\documentclass{amsart} 
\usepackage{latexsym,amssymb,graphicx,amscd}
       
\newtheorem{thm}{Theorem}
\newtheorem{lem}[thm]{Lemma}
\newtheorem{prop}[thm]{Proposition}

\newtheorem{ex}[thm]{Example}
\newtheorem{que}[thm]{Question}

\newcommand{\BZ}{{\mathbb{Z}}}
\newcommand{\BQ}{{\mathbb{Q}}}

\newcommand{\BO}{{\mathcal{O}}}

\DeclareMathOperator{\Tor}{Tor}

\DeclareMathOperator{\kernel}{kernel}

\DeclareMathOperator{\Int}{Int}

\begin{document}
\title{Quantum invariants and free $\BZ_{p^2}$--actions on 3--manifolds}.

\author{ Patrick M. Gilmer}
\address{Department of Mathematics\\
Louisiana State University\\
Baton Rouge, LA 70803\\
USA}
\email{gilmer@math.lsu.edu}
\thanks{The first author was partially supported by NSF--DMS--0604580}
\urladdr{www.math.lsu.edu/\textasciitilde gilmer/}

\author{ Khaled Qazaqzeh}
\address{Department of Mathematics, Yarmouk University, Irbid, 21163,
Jordan} \email{qazaqzeh@math.lsu.edu}
\urladdr{http://faculty.yu.edu.jo/qazaqzeh/}

\date{June 16, 2008}

\begin{abstract}  We give a congruence for the quantum invariant of a $\BZ_p$--quotient of  a $3$--manifold with a $\BZ_{p^2}$ action.  
We show the congruence does not hold for quotients of 3--manifolds with a  $\BZ_{5} \times \BZ_{5} $ action.  
\end{abstract}

\maketitle

\section{Introduction}

Let $p$ denote an odd prime. 
The most simple kind of finite cyclic covers are those which are quotients of infinite cyclic covering spaces. 
We call such covers { \it simple} cyclic covers.  Such covers are formed by stacking slit copies of the base. As  TQFTs satisfy nice properties with respect to stacking, one can calculate quantum invariants of  finite simple cyclic covers  nicely from data for the base and the covering \cite{G1}.  Moreover one obtains in this way congruences modulo $p$  for the quantum invariants of  simple $\BZ_p$--cyclic covers of closed oriented connected 3--manifolds \cite{G2,G3}: 

\begin{thm} \label{infinite} Suppose the infinite  group $\BZ$  acts freely and preserving the orientation on  a connected oriented 3--manifold $\check N$ with a  compact quotient. Then there exist  $m$ and $n$  in $\BZ$  with
\[ \left\langle  \check N /p \BZ \right\rangle_p\equiv {\kappa}^m n \pmod{p\  \BO_p}.\]
\end{thm} 

 Here
$\BO_p$ denotes $\BZ[A,\kappa]$ where $A$ is a primitive $2p$ th root of unity and $\kappa^2= A^{-6-p(p+1)/2},$ and
  $\left\langle \  \right\rangle_p \in  \BO_p[1/p]$ denotes  the invariant   \cite{BHMV} of closed oriented 3--manifolds possessing an integral weight, as well  a (possibly empty) $p$--admissibly colored (with colors integers in range $[0,p-2]$) fat graph. 
 This invariant  can be computed using the TQFT $(Z_p,V_p)$  of  \cite{BHMV}. We use a modified form of this TQFT  as in \cite{G2} where $p_1$--structures are replaced by integral weights on 3--manifolds and Lagrangian subspaces of the first homology of surfaces. When one raises this weight, one multiplies the invariant of a closed 3--manifold by the phase factor  $\kappa$.  One has that $\left\langle S^1 \times S^2 \text{ with weight zero} \right\rangle_p =1$, but $\left\langle S^3  \text{ with weight zero}  \right\rangle_p \notin \BO_p$.   
 Thus this is a different normalization for the invariant than used in \cite{M}.

It is a natural to wonder whether Theorem \ref{infinite} continues  to hold for a more general class of  $p$--fold covers $\tilde N$ of a closed oriented 3--manifold $N$. However this congruence would only be possible if the values of the quantum invariant  of  $\tilde N$ lies at least in  $\BO_p.$

\begin{prop} Assume $p \ne 3.$  Let $\tilde N$ be a $p$--fold cyclic cover of a closed oriented 3--manifold $N$ with empty  colored fat graph.  The following are equivalent.
\begin{itemize} 
\item[(1)] $\left\langle\tilde N  \right\rangle_p \in  \BO_p$
\item[(2)] $H^1(\tilde N, \BZ_p) \ne 0$
\item[(3)] There is a $p^2$--fold regular cover $\check N$ of $N$ with group of covering transformations $G$ isomorphic to  either $\BZ_{p^2}$  or  $\BZ_{p} \times \BZ_{p}$ such that 
$\tilde N$ is the quotient of $\check N$ by a subgroup of $G$ of order $p.$
\end{itemize}
\end{prop}

\begin{proof}
According to \cite{M}, $\left\langle\tilde N  \right\rangle_{p} \notin  \BO_p$ if $M$ is a $\BZ_p$--homology sphere. Thus (1) implies (2).
By  \cite{CM,G4} (2) implies (1). 

Clearly (3) implies (2). For the converse,  let $\psi$  be an eigenvector for the action of a generator of the group of covering transformation $\tilde N \rightarrow N$ on $H^1( \tilde N, \BZ_p)$.
 As the order of this generator is $p$, then any eigenvalue is a $p^{th}$ root of unity in $\BZ_p$.
 But the only element with order a divisor of $p$ in $\BZ_p^*$ is one. Therefore, one is the only  eigenvalue for  this generator.  Thus 
  $\psi$ is fixed by the group of covering transformations. Let $ \check N$ be the $\BZ_{p}$ cover of $\tilde N$ which is classified by $\psi$. The cover $\check N \rightarrow N$ is then regular with group of covering transformations a group of order $p^2$. As is well-known, a group of order $p^2$ is either $\BZ_{p^2}$, or $ \BZ_{p} \times \BZ_{p}$. 
\end{proof}

All the implications in the above proof hold also in the case $N$ contains a non-empty colored fat graph  except perhaps (1) implies (2). Similarly in the case $p=3$, all these results hold except (1) implies (2). We remark that the results of this paper which concern quantum invariants are uninteresting in the case $p=3,$ as $\left\langle N \right\rangle_3$ is always some power of 
$\kappa$. However some of the purely topological results may be of interest in the case $p=3$.

Thus we are lead to:

\begin{que}\label{ques} Suppose  $G$,
a group of order $p^2$,  acts freely on  a closed connected  oriented 3--manifold $\check N$. Let $H\subset G$ denote a subgroup of order $p$. Must there exist  $m$ and $n$  in $\BZ$  with \[ \left\langle \check N /H \right\rangle_p\equiv {\kappa}^m n \pmod{p\  \BO_p} ? \]
\end{que} 

We will show that the answer is ``yes'' if $G \approx  \BZ_{p^2}$, and ``no'' if $G \approx  \BZ_{5} \times \BZ_{5}$.

The example that we give for $\BZ_{5} \times \BZ_{5}$ when modified in the most natural way will not show the answer is ``no'' for any $\BZ_{p} \times \BZ_{p}$  with $p>5$.

The Appendix  discusses a case which was omitted in the proof of a related congruence in previous paper of the first author \cite{G3}.

\section{Simple covers and weak type--$p$  surgery  \label{bochs}}
Before the proof of Theorem \ref{main} in the next section, we must prepare some tools.
We let $X$ denote the base of a $\BZ_k$--cyclic cover where $k$ is a power of $p$.   Then we have
$\chi \in H^1( X, \BZ_k)$ which classifies the cover together with a choice of generating covering transformation. Consider  the Bockstein homomorphism 
 $\beta_k$ associated to the short exact sequence of coefficients:
\[ 0 \rightarrow \BZ \rightarrow \BZ \rightarrow \BZ_k \rightarrow 0.\]
It fits into a long exact sequence.
\[ \begin{CD}  H^1( X) @>{\times k}>> H^1( X)@>{}>> H^1( X, \BZ_k) @>{\beta_k}>> H^2( X) \end{CD}.\]
 Clearly $\chi \in H^1( X, \BZ_k)$ classifies a simple cover if and only if $\beta_k(\chi)=0.$ 
 
 { \it For the rest of this paper $N$ will denote a closed oriented connected $3$--manifold}. 
Recall that  the Poincare duality  isomorphism from $H^2(N)$ to $H_1(N)$ is given by capping with the fundamental class: $\cap [N]$. Let $\rho$ denote the inverse isomorphism from $H_1(N)$ to $H^2(N).$

\begin{lem}\label{Gamma1} Suppose $\mathcal{G} = \cup_i^n \gamma_i$  is a link.      The cover given by $\chi$ restricted to the complement of $\mathcal{G}$ is simple if and only if $\beta_{k}(\chi)$  is in the span of $\{ \rho([\gamma_i] \}$. The cover restricted to $\gamma_i$ is characterized by $\chi(\gamma_i) \in \BZ_k$.
\end{lem}

\begin{proof}  Let  $\nu(\mathcal{G})$ denote a closed tubular neighborhood of $\mathcal{G}$. Let  $N_{\mathcal{G}}$ denote $N$ with the interior of $\nu(\mathcal{G})$ deleted. 
Consider the long exact sequence of the pair $(N,N_{\mathcal{G}})$:
\[ \begin{CD}  H^2( N,N_\mathcal{G}) @>j>> H^2( N)@>{}>> H^2(  N_\mathcal{G})   \\
@A\text{excision isomorphism inverse}AA \\
H^2(\nu(\mathcal{G}), \partial \nu(\mathcal{G})) @.\hspace{-30pt}\nearrow \hspace{30pt}\\
@A\text{Thom  isomorphism}AA \\
H^0(\mathcal{G})
 \end{CD}.\] 
 We  have that  $H^0(\mathcal{G})$ is free on the components. 
 The diagonal map is defined so that
the diagram commutes. Using the  dual chain complex construction of the Poincare duality isomorphism e.g. \cite[Thm 65.1]{Mu}, we see that the  image of the generator of $H^0(\gamma_i)$ under  the diagonal map is
  $\rho[\gamma_i]$. 
Thus the image of $j$ is spanned by the $\rho [\gamma_i]$ in $H^2(N)$. As the horizontal sequence is exact, 
$\beta_{k}(\chi)$ restricted to $H^2(N_\mathcal{G})$ is zero.  The cover of the simple closed curve is classified by the map
 from the infinite cyclic group $H_1(\gamma)$ to  $\BZ_k,$ which maps the generator to $\chi(\gamma_i)$.
\end{proof}

\begin{ex} \label{exL} A connected $\BZ_p$--cover of the lens space  $L(p^2,1)$ is simple on the complement of a simple closed curve representing $p$ times a generator of $H_1(L(p^2,1))$, 
but the covering of the curve is trivial. If instead we consider a simple closed curve representing a generator of $H_1(L(p^2,1))$, this same covering is simple on the complement of the curve and non-trivial on  the curve. However, if we consider connected sum of $p$ copies of 
$L(p^2,1)$ and the cover given by  a character $\chi$ which sends the standard generator on each summand to $1 \in \BZ_p$, then there is no simple closed curve which is covered non-trivially but whose complement is covered simply. 
\end{ex}

We now wish to learn how to pick ${\mathcal{G}} = \cup_i^n \gamma_i$  as above with control over $\chi(\gamma_i)$.   To understand the general case,  we must discuss the relation between the Bockstein and the linking form. We must also relate our characters in 
$H^1(N, \BZ_k)$ to characters in $H^1(N, \BQ/\BZ).$   For this last issue, we have a commutative diagram of coefficients with short exact rows.

\[ \begin{CD}  0@>>>\BZ @>{\times k}>> \BZ@>>> \BZ_k @>>> 0\\
			@.     	@|  @V{\times \frac 1 k}VV @VVV @.    \\
		0@>>>\BZ @>>> \BQ@>>> \BQ/\BZ @>>> 0
		\end{CD}
\]

The residue class of one in $\BZ_k$ is sent to $1/k \in \BQ/\BZ$. We let the Bockstein associated to the lower sequence be denoted $\beta$.  We obtain in this way a commutative diagram with exact rows: 

\[ \begin{CD}   H^1( N, \BZ_k) @>{\beta_k}>> H^2( N)  @>{\times k}>> H^2( N)\\
@VVV @|  @VVV\\
 H^1( N, \BQ/\BZ) @>{\beta}>> H^2( N)  @>{}>> H^2( N,\BQ)
\end{CD}. \]

Thus given a character $\chi: H_1(N) \rightarrow \BZ_k$, if we compose with the standard inclusion 
$ \BZ_k$
to $\BQ/\BZ$ and then apply $\beta$,  we get $\beta_k(\chi)$. Also $\beta$ maps onto the torsion subgroup of $H^2( N)$, denoted $\Tor (H^2(N)) $. Note that $ \beta ( \chi ) \cap[N]$ is the Poincare dual of the Bockstein of $\chi.$
We are interested in the bilinear form:
\[ b: H^1(N,\BQ/\BZ) \times H^1(N,\BQ/\BZ) \rightarrow \BQ/\BZ\]
given by
\[b( \chi_1, \chi_2) = \chi_1(  \beta ( \chi_2) \cap[N]).\]
So using \cite[p.254]{S}
\[b( \chi_1, \chi_2) =  \chi_1\cap (\beta ( \chi_2 ) \cap[N]) =  (\chi_1\cup \beta ( \chi_2 )) \cap[N].\]

\begin{lem}  The form  $b$  is  symmetric. 
\end{lem}
\begin{proof} 
 \begin{align*} 
0 &= \beta (\chi_1 \cup \chi_2) \text{ by exactness as $H^3(N,\BZ) \rightarrow H^3(N,\BQ) $ is  one to one}\\
 & =   \beta ( \chi_1 ) \cup \chi_2 -\chi_1  \cup \beta ( \chi_2 ) \text{ as $\beta$ is a derivation}\\
 &=   \chi_2 \cup\beta ( \chi_1 )  -\chi_1  \cup \beta ( \chi_2 ) \text{.  So } \end{align*}
 \[
 b(\chi_2, \chi_1) -  b(\chi_1, \chi_2) =   (\chi_2 \cup  \beta ( \chi_1)) \cap [N]-(  \chi_1 \cup \beta ( \chi_2 )) \cap [N] =0   \]
\end{proof}

As $b$ vanishes on $\kernel(\beta) \times H^1(N,\BQ/\BZ)$ and  $H^1(N,\BQ/\BZ)/ \kernel(\beta) \approx \Tor (H^2(N))$, there is an induced  nonsingular symmetric form  \cite[p.75--76]{DK}
\[  {\mathfrak b}: \Tor (H^2(N)) \times \Tor (H^2(N)) \rightarrow \BQ/\BZ\]
such that 
\[ b( \chi_1, \chi_2) =  {\mathfrak b}(\beta(\chi_1), \beta(\chi_2)) .\] In fact, 
 transferred to $ \Tor (H_1(N))$ by the Poincare duality isomorphism, one can check that ${\mathfrak b}$  becomes the well known linking form on $\Tor (H_1(N)).$

Wall \cite{W} has classified non-singular symmetric  forms on finite abelian groups. Such forms have a orthogonal decomposition into primary summands associated to each prime. The $p$--primary summand is isomorphic to a direct sum of elementary forms of  two types: $A_{p^t}$ and $B_{p^t}$. The underlying group of  both is $\BZ_{p^t}$. The form on  $A_{p^t}$ sends $(x,y)$ to $xy/p^t$. The form on  $B_{p^t}$ sends $(x,y)$ to $n_t xy/{p^t}$, where $n_t$ is a non-square unit in $\BZ_{p^t}.$  It is pleasant that we do not have  to deal with the even prime, which is also studied by Wall but which is more complicated.

\begin{lem}\label{scc} If $\chi \in H^1(N,\BZ_{p})$ is nonzero, then we can pick a link ${\mathcal{G}} = \cup_i^n \gamma_i$   in $N$, so that  the cover given by $\chi$ restricted to the complement of $\mathcal{G}$ is simple and $\chi$ on each $\gamma_i$ is nonzero.
\end{lem}

\begin{proof}   We reinterpret $\chi$  to lie in  $H^1(N,\BQ/\BZ)$.  We can view $\beta(\chi)$  as  an element of the $p$--primary subgroup of $\Tor (H^2(N))$ which we assume is already decomposed and identified in the way described by Wall. Let $n$ be the number of summands where $\beta(\chi)$ projected into that summand is nonzero. The image of $\beta(\chi)$ under each of these projections has order  $p$.  It is always  possible to find an  element $x_i$ of this summand  which pairs  under the torsion form with the projection of $\beta(\chi)$ to yield  $1/p \in \BQ/\BZ$. Let $\gamma_i$ be Poincare dual to $x_i$.
Then  $\beta ( \chi ) $ is in the span of the
Poincare duals of the $\{[ \gamma_i]\}$. Moreover for each $1 \le i \le n$,  
\[
\chi(\gamma_i)= b( \chi, \beta^{-1}(x_i) )={\mathfrak b} (\beta(\chi) , x_i)=1/p. \]
\end{proof}

If $\gamma$ is a simple closed curve  in $N$, and $\mathcal{T}_{\gamma}$ is the boundary of a tubular neighborhood $\nu_{\gamma}$ of $\mathcal{T}_{\gamma}$, let $\mu_{\gamma} \subset \mathcal{T}$ denote a meridian of $\gamma$. We may also pick a longitude $\lambda_{\gamma}$. This is a curve on $\mathcal{T}_{\gamma}$ which is homologous to $\gamma$ (with some orientation) in  $\nu_{\gamma}.$  
For $d$ an integer bigger than 1, by weak type--$d$ surgery  along  $\gamma$ in $N$ \cite{G4},  we mean the process of removing  $\nu_{\gamma}$ and regluing  it so a curve on $\mathcal{T}$ homologous to $n \mu +d l \lambda$ for some integers $n$, and $l$ bounds in the reglued solid tori. 
Suppose $\lambda'$ is another choice of longitude. Then $\lambda'= \lambda +x \mu$, so $n \mu +d l \lambda=(n- d l x) \mu +d l \lambda'$.
Thus 
the notion of weak type--$d$ surgery does not depend on which longitude was chosen. A weak type--$d$ 
surgery can be undone with another. The equivalence relation generated by weak  type--$d$ surgery is called weak $d$--congruence.

The following Proposition  is a weakened form of \cite[Prop 2.14]{G4}. We include it here as the proof  helps 
to prove  Proposition \ref{cut} and to motivate the proof of  Lemma \ref{p2equiv}. 

\begin{prop}\label{pequiv} If  $H_1(N,  \BZ_p)$ is non-zero,  then  $N$ is weakly $p$--congruent to a 3--manifold $M$ with $H^1(M )$ non-zero.
 \end{prop}

\begin{proof} We  pick a nonzero  $\chi \in H^1(N,\BZ_{p})$ and apply
Lemma \ref{scc}. Let $\hat \chi \in H^1(N\setminus \mathcal{G},\BZ)$ be a lift  to $\BZ$ of the restriction of $\chi$ to $N\setminus \mathcal{G}.$ Let $ T_i$ be the boundary of a tubular neighborhood of $U_i$ of $\gamma_i$ equipped with a meridian $\mu_i$, and   a chosen longitude  $\lambda_i \in H_1(T_i)$. The character $\hat \chi$ induces a character $\hat \chi_i:H_1(T_i) \rightarrow \BZ$. Suppose the image of $\hat \chi_i$ is $k_i \BZ$. Note that $k_i$ and $p$ must be relatively prime. Then 
$\hat  \chi_i([\mu_i])= p s_i k_i$ and $\hat  \chi_i([\lambda_i])= -n_i k_i$ 
for integers $s_i,k_i, n_i$ with  $s_i, n_i$ relatively prime  and with $n_i$ prime to $p.$ Then $\hat  \chi_i$ vanishes on the homology class of a curve $\mu'_i$ representing $n_i [\mu_i]+ p s_i [\lambda_i]$. Thus $\hat  \chi$ extends uniquely to $M$ obtained by doing a weak type--$p$ surgery along each curve $\gamma_i$. As this extension will be nonzero, $H^1(M,\BZ)$ is nonzero. By definition, $M$ is weakly $p$--congruent to $N$.
\end{proof}

By a $p$--surface $F$, we mean \cite{G4} the result of attaching, by a map $q$, the whole
boundary of an oriented surface $\hat F$ to a collection of circles $\{ S_i\}$ by a map  which when restricted to the  inverse image under $q$ of each $S_i$ is a $p  t_i $--fold ( possibly disconnected)  covering space of $S_i.$  If each component of each $q^{-1} S_i$ is itself a covering space of $S_i$ with degree divisible by $p$, we say $F$ is a good $p$--surface. The $p$--cut number of a 3--manifold $N$, denoted $c_p(N)$, is the maximum number of disjoint piecewise linearly embedded good $p$--surfaces that we can place in $N$ with a connected complement. 
The following Proposition allows us to interpret \cite[Theorem 4.1]{G4} as a corollary of \cite[Theorem 4.2]{G4}.

\begin{prop}\label{cut} If  $H_1(N,  \BZ_p)$ is non-zero,  then   
 $c_p(N)>0.$ 
 \end{prop}
\begin{proof} Continuing with the proof of Proposition \ref{pequiv},  let $\chi'= (1/k) \chi$,  $U=\cup U_i$ appearing in this proof.   In  $N\setminus \Int U$, we can find a connected surface $F$ with a connected complement which is dual to $\chi'_{|U}$ and meets  each $T_i$ in $\mu'_i$. $F$ maybe completed to a good $p$--surface with
a connected complement by adding the mapping cones of the projections from $\mu'_i$ to $\gamma$
\end{proof}

\section{$G= \BZ_{p^2}$}

\begin{thm} \label{main} Suppose $\BZ_{p^2}$ acts freely on  closed oriented   connected  3--manifold $\check N$. If $H$ is the subgroup of   $\BZ_{p^2}$ of order $p$,  then  there exists  $m$ and $n$  in $\BZ$  with
\[ \left\langle\check N/H  \right\rangle_p\equiv {\kappa}^m n \pmod{p\  \BO_p}  \]
\end{thm}

\begin{proof} In this situation, we say $ \check N/H$ is a somewhat simple $\BZ_p$--cyclic cover of $\check N/\BZ_{p^2}$, which we will denote by $N$. A basic example of somewhat simple $\BZ_p$--cyclic cover which is not a simple $\BZ_p$--cyclic cover is the lens space $L(p,q)$ which is a somewhat simple $\BZ_p$--cyclic cover of $L(p^2,q)$.

  We need a version of  Lemma \ref{scc} for characters in $H^1(N,\BZ_{p^2})$. 

\begin{lem}\label{scc2} If $\chi :H_1(N) \rightarrow \BZ_{p^2}$ is an epimorphism, then we can pick a link ${G} = \cup_i^{n} \gamma_i $   in $N$, so that  the cover given by $\chi$ restricted to the complement of $G$ is simple and for each $i$, $\chi(\gamma_i)$ is either one, or is  $p$.
\end{lem}

\begin{proof}  Again we reinterpret $\chi$  to lie in  $H^1(N,\BQ/\BZ)$, view $\beta(\chi)$  as  an element of the p--primary subgroup of $\Tor (H^2(N))$ which we assume is already decomposed and identified in the way described by Wall.

 Let $0 \le i \le n$ index the summands where $\beta(\chi)$ projected into that summand is  nonzero. Let $\beta_i$ the projection of $\beta(\chi)$ in the $i$th summand.  
 If $\beta_i$ has order $p^2$,  it is  possible to find a  element $x_i$ of this summand  which pairs  under the torsion form with $\beta_i$ to yield   $1/p^2 \in \BQ/\BZ.$
 If $\beta_i$ has order $p$, we can find a  element $x_i$ of this summand  which pairs  under the torsion form with $\beta_i$ to yield   $1/p \in \BQ/\BZ.$ Let $\gamma_i$ be Poincare dual to $x_i$.
Then  $\beta ( \chi ) $ is in the span of $\{\rho[ \gamma_i]\} $. Moreover for each $1 \le i \le n$,  $\chi(\gamma_i)$ is either $1/p^2$ or $1/p$.  Now reinterpret $\chi$ to take values in $\BZ_{p^2}$.\end{proof}

\begin{lem}\label{p2equiv} Suppose $\chi :H_1(N) \rightarrow \BZ_{p^2}$  is  an epimorphism.  Let $\tilde N$ denote the associated $p$--fold cover classified by $\pi_p \circ \chi$, where $\pi_p: \BZ_{p^2} \rightarrow \BZ_p$ is reduction modulo $p$.
Then $N$ is weakly $p$--congruent to a manifold  $M$ with a simple 
$\BZ_p$ covering $\tilde M$.  $\tilde M$ may be obtained by 
 weak type--$p$ surgery on a link in $ \tilde N$.
\end{lem}   
 
\begin{proof} Apply
Lemma \ref{scc2}. Let $\epsilon_i \in \BZ$ be $1$ if $\chi(\gamma_i)=1$  and  $\epsilon_i$ be $p \in \BZ$ if $\chi(\gamma_i)=p$.
Let $\hat \chi \in H^1(N\setminus {\mathcal{G}},\BZ)$ be a lift of $\chi$  to $\BZ.$    Then 
$\hat  \chi_i([\mu_i])= p^2 s_i k_i$ and $\hat  \chi_i([\lambda_i])= - \epsilon_i n_i k_i$ 
for integers $s_i,k_i, n_i$ with  $s_i, n_i$ relatively prime  and with both $n_i$ and $k_i$ prime to $p.$ Then $\hat  \chi$ vanishes on the homology class of a curve $\mu'_i$ representing $n_i [\mu_i]+ 
\frac {p^2} {\epsilon_i} s_i [\lambda_i]$.

Thus we can do weak $\frac {p^2} {\epsilon_i}$--type surgery to $N$ along each $\gamma_i$ to obtain a weakly $p$--congruent manifold $M$ such that  $\hat \chi$ extends to $\hat \chi': H_1(M) \rightarrow  \BZ.$ This induces a simple 
$\BZ_p$--cover $\tilde M$ of $M.$ 
We have that $\tilde M$  is obtained from $\tilde N$ by a sequence of surgeries along the curves in $\tilde N$ which lie over the $\gamma_i.$
 If  $\epsilon_i=1$ then $\gamma_i$ is  covered by a single  curve  in $\tilde N$ and we  perform weak type--$p$ surgery along  this curve.
 If 
$\epsilon_i=p$ then $\gamma_i$ is  covered by $p$ disjoint  curves  in $\tilde N$ and we  perform weak  type--$p^2$  surgeries ( which are also weak type--$p$ surgeries ) along each of these curves. \end{proof}

  By \cite[Theorem 3.8]{G4},  $\left\langle \tilde N \right\rangle_p$ is, up to phase,  \[\left\langle \tilde M \text{ with some $\BZ_p$--equivariant integrally colored framed link } \right\rangle_p.\]  Since $\tilde M$ is a simple $\BZ_p$--cover,  by \cite[11.1]{G2}, 
this last expression  must satisfy the stated congruence.
\end{proof}

\section{$G= \mathbb{Z}_{p} \times \mathbb{Z}_{p}$}

Let $L_p$ be the two component link obtained by replacing one component of a Hopf link with a $(p,1)$ cable with framing $p$ on the cabled component and framing zero on the 
uncabled component.  Let $M_p$ be obtained  by performing framed surgery on $L_p$.
The linking matrix of the framed link  is 
$\begin{bmatrix}
0 & p\\
p& p
\end{bmatrix}$.
Thus $H_1(M_p) = \BZ_5 \oplus \BZ_5$. 
Let $\check M_p$ denote the maximal abelian cover of $M$, i.e. the covering space whose fundamental group is the commutator subgroup. 
It is a regular $\BZ_p \oplus \BZ_p$ cover of $M$. Let $L_p'$ be a $(p, p)$--torus link  with an unknot going around it (a satellite of the Hopf link with the (5, 5) torus link the  companion of one component of the Hopf link.). We give $L'$ the framing +1 on each component of the cabled part and the zero framing  on the uncabled component.  Let $\tilde M_p$ be the 3--manifold obtained from $S^{3}$ by doing surgery along $L'_p$ with weight zero. $\tilde M_p$ is a  $\BZ_p$--cyclic cover of $M$, and thus a quotient of $\check M$ by an order $p$ cyclic subgroup of $\BZ_p \oplus \BZ_p$.

Below will  use the $z^i$ basis \cite{BHMV1} for the Kauffman skein of a genus one handlebody.

Let $H_n$ denote the Kauffman bracket of an $n$ component positive Hopf link where every component has framing one. By \cite[Prop 11 (ii)]{L2}, we have for $n\ge 1$,
\[H_n= (A^2-A^{-2})^{-1}
\sum_{r=0}^{n-1} \binom {n-1}{r} A^{(n-2r+1)^2-1} (A^{2(n-2r+1)}-A^{-2(n-2r+1)})\]
Note $H_0=1$, so this is not valid for $n=0.$ 

\begin{thm} \label{ex}  $\BZ_5 \times \BZ_5$ acts freely on  $\check M_5$. There is an order $5$ cyclic subgroup $H$ of $\BZ_5 \times \BZ_5$  such that there are no $m$ and $n$  in $\BZ$  with
\[  \left\langle \check M_5/H  \right\rangle_5 \equiv  {\kappa}^m n \pmod{5\ \BO_5}  \]
\end{thm}

\begin{proof}
We let $\BO_5= \BZ[\zeta_{20}]$ with $A= \zeta_{20}^2$ and $\kappa= \zeta_{20}^{-1}. $ Then  by formula for $\eta$ near the beginning of \cite[\S 2]{BHMV}, 
$\eta= \frac 1{5}( {2 \zeta_{20} + \zeta_{20}^3 + \zeta_{20}^5-3\zeta_{20}^7}).$ 

One has $\Omega_{5} = 1 +\delta z$, where $\delta  =-A^{-2}-A^2$. Replacing a component with framing one with
$t^{-1}\Omega_{5} = 1 -A^{-3}\delta z$, has the same effect as  first changing the framing to zero and  then replacing with 
 $\Omega_{5}$. Here $t$ denotes the twist map on the Kauffman  skein of 
 of the solid torus \cite{BHMV1}. 
We let $\left\langle L'(\Omega_{5})\right\rangle$ be the Kauffman bracket evaluation  of the linear combination over $\BO_5$ that we obtain if we replace each component of $L'$ by $\Omega_{5}$. One has that $\left\langle \tilde M_5 \right\rangle_5 =\eta_{5}^{7} \left\langle L'(\Omega_{5})\right\rangle$.  

To compute $\left\langle L'(\Omega_{5})\right\rangle$, we compute a linear combination of $2^6$ brackets of Hopf links with zero to six components and varying framings. If we expand out the zero framed component first, we obtain: 
\[\left\langle L_5'(\Omega_{5})\right\rangle
=\sum_{k=0}^5 \binom 5 k \delta^k H(k) - A^{-3} \delta \sum_{k=0}^5 \binom 5 k \delta^k H(k+1)\]

We obtain
\[\left\langle \tilde M_5 \right\rangle_5 =  -2\zeta_{20} +4\zeta_{20}^{3} - \zeta_{20}^{5} - 2\zeta_{20}^{7} \]
Comparing this with the finite  list of elements of $\BO_5/5 \BO_5$ which are images of numbers of the form
$n \kappa^{m}$
for any $n,m\in \mathbb{Z}$, we conclude that  $\left\langle \tilde M \right\rangle_5$ does not have this form.\end{proof}

Now we consider $p=7$. We let $\BO_7= \BZ[\zeta_{14}]$ with $A= \zeta_{14}$ and $\kappa= A^4. $ Then  
$\eta= \frac 1{7}( {-2 \zeta_{14}- \zeta_{14}^2-2 \zeta_{14}^3 + 2\zeta_{14}^4 +\zeta_{14}^5}).$ 
We have 
$\Omega_{7} = (2-\delta^2)  +\delta z + (\delta^2-1) z^2$, and $t^{-1} \Omega_{7} = (1+ A^{6}- A^{6}\delta^2)  - A^{11}\delta z +A^{6} (\delta^2-1) z^2.$ We have that
$\left\langle \tilde M_7 \right\rangle_7 =\eta^{9} \left\langle L_7'(\Omega_{7})\right\rangle$.  Using multinomial coefficients, we have 

\begin{align*} \left\langle L_7'(\Omega_{7})\right\rangle
=&  (1+ A^{6}- A^{6}\delta^2) \sum_{i+j+k=7} \left(\frac 7{i,j,k}\right)  (2-\delta^2)^i \delta^j (\delta^2-1)^k H(j+2k) \\
& - A^{11}\delta \sum_{i+j+k=7} \left(\frac 7{i,j,k}\right)  (2-\delta^2)^i \delta^j (\delta^2-1)^k H(j+2k+1)\\
& +A^{6} (\delta^2-1) \sum_{i+j+k=7} \left(\frac 7{i,j,k}\right)  (2-\delta^2)^i \delta^j (\delta^2-1)^k H(j+2k+2).
\end{align*}

It follows that:
\[ \left\langle \tilde M_7 \right\rangle_7= 7  (176993 + 397520  A - 318640  A^2 - 220548  A^3 - 
        98084  A^4 + 495621 A^5) \]

Thus $\left\langle \tilde M_7 \right\rangle_7 \in 7  \BO_7$.       
Notice that $H_1(\tilde M_p)= \BZ^{p-1}.$ By \cite[4.3]{CM}
\[ \left\langle \tilde M_p \right\rangle_p \in  (1-\zeta_p)^{\lceil(p-1)(p-3)/6-(p-3)/2\rceil} \BO =
(1-\zeta_p)^{\lceil(p^2-7p +12)/ 6 \rceil} \BO. \]
As $\lceil(p^2-7p +12)/ 6 \rceil \ge p-1$ for $p \ge 11$, we have  $\left\langle \tilde M_p \right\rangle_p \in p  \BO_p$, for $p \ge 11$. Thus these examples for $p \ge 7$ are consistent with a ``Yes'' answer to Question \ref{ques}. To find  examples allowing one to say ``No'', one should begin by looking for  3--manifolds with first homology $\BZ_p \times \BZ_p$ with a p--fold cover with the dimension of first homology with $\BZ_p$--coeffients less than
$6 \frac {p-1}{p-3} +3.$ Of course $\tilde M_7$ did satisfy this
equation but did not provide a ``No'' answer.

We used Mathematica \cite{Wo} for many of the computations in this section.

\appendix

\section{Covers of degree prime to the order of the phase factor by Patrick M. Gilmer}

In the proof of \cite[Theorem 1$'$]{G2},  on page 171 Proposition 2 only applies when the cover restricted to $\gamma$ is non-trivial.   The case where the cover restricted to $\gamma$ is trivial was not addressed.  In view of Example \ref{exL}, this case may arise. 

Theorem 1$'$  concerns the quantum invariants $\left\langle N \right\rangle_{2r}$ where $r$ is prime to  $p$. We continue to assume $p$ is an odd prime. The congruence given is not modulo phase but exact. As the extra structures used to resolve the ``framing anomaly ''   in \cite{G2} are the $p_1$--structures of \cite{BHMV} rather than integral weights and Lagrangians,  3--manifolds, in this appendix, will be equipped with $p_1$--structures. We note that $p_1$--structures are more natural in this context, as covering spaces are equipped with $p_1$--structures induced from the base. Also 3--manifolds are allowed to have possibly empty admissibly colored fat graphs.  The quantum invariant $\left\langle \  \right\rangle_{2r}$ takes values in $\BZ[ \frac 1 {2r} ,\xi]$ where 
$\xi$ is a primitive $4rth$ root of unity if $r$ is even and a primitive  $8rth$ root of unity if $r$ is odd. The colors of this theory are from the set $C$ of  the integers from zero to $r-2$.  The following result addresses the missing case in the proof of  \cite[Theorem 1$'$]{G2}. In fact  \cite[Theorem 1$'$] {G2} has the same conclusion as Theorem \ref{app} under the weaker hypothesis:
$N_p$ is a connected $\BZ_p$--covering of $N$.

\begin{thm} \label{app} Let  $N_p$ be a connected non-simple $\BZ_p$--covering of $N$ given by $\chi \in H^1(N, \BZ_p)$. 
Suppose $\gamma$ is a simple closed curve in $N$ such that $\chi[\gamma]=0,$ and the cover restricted 
to $N \setminus \gamma$ is simple, then 
\[ \left\langle  N_p \right\rangle_{2r} \equiv  \left\langle  N \right\rangle_{2r}^p \pmod{ p \ \BZ[ \frac 1 {2r} ,\xi]} \]
\end{thm}
\begin{proof} Let  $\nu$ denote a closed tubular neighborhood of $\gamma$, ${\mathcal T} = \partial \nu$, $E =-N \setminus  \Int(\nu)$. Let $E_p$ be  the given cover  restricted to $E$. This cover restricted to  
${\mathcal T}$ is a disjoint union of $p$ tori
 $\mathcal T_p=\coprod_{i=1}^p {\mathcal T}_i.$  
 Similarly this cover restricted to  
${ \nu}$ is a disjoint union of $p$ solid tori
 $\nu_p=\coprod_{i=1}^p \nu_i$ with $\partial \nu_i= \mathcal T_i.$
 We index these so the covering  transformation specified by $\chi$, sends $\nu_i$ to $\nu_{i+1}$ for all $i$.
 
Let $\hat \chi: H_1(E)  \rightarrow \BZ$ be a surjective lift of $\chi.$ 
Let $\mu'$ be a simple closed curve in $\mathcal T$ which generates the kernel of  $\hat \chi$ composed with the map induced by  inclusion $H_1(\mathcal T) \rightarrow H_1( E)$.  Let $H$ be a  handlebody  with boundary  $\mathcal T$ such that $\mu'$ bounds a disk in $H$.  We let $N'= E \cup_{\mathcal T} H$. The cover on $E$ extends to a simple cover $N'_p$ of $N'$.
Let $H(j)$ denote $H$ with  the core colored $j$,
$N'(j)=  E \cup_{\mathcal T} H(j)$, 
and $N'(j)_p$ denote $N'_p$  with the $p$ circles covering the core all colored 
$j$.  We let $H_i$ denote the component of the  cover of $H$ with $\partial H_i= \mathcal T_i,$ and 
let  $H(j)_i$ denote $H_i$ with the core colored $j$.
We will use $\{ [H(j)] | 0 \le j \le  r-1 \}$ as a basis for $V_p(\mathcal T)$.
It is orthogonal with respect to the Hermitian form $\left\langle \ , \  \right\rangle_\mathcal T$   
on $V_p(\mathcal T)$.

Let $\mathcal S$ denote the set of sequences  of colors of length $p$. We denote the $i$th term of the  sequence $\sigma$ by  $\sigma_i$.  Let $\tau: \mathcal S \rightarrow \mathcal S$ be the cyclical shift map with 
$\tau(\sigma)_i=  \sigma_{i-1\pmod{p}} $.
The orbits of $ \mathcal S$ under powers of the shift map are of two types. There are singletons
given by constant sequences, and orbits with cardinality $p$ made up of non-constant sequences.
We index the constant sequences by the set of colors $C$. We denote the sequence which is constantly $j$ by $\tilde j$.

 If $\sigma \in \mathcal S$, let $H(\sigma)$ denote
\[H(\sigma_1)_1 \otimes H(\sigma_2)_2 \otimes \cdots H(\sigma_p)_p \in V_p(\coprod_{i=1}^p \mathcal T_i)= V_p(\mathcal {T}_p).\]

Using orthonormality, 
\begin{equation}\label{1} [\nu]=  \sum_{j \in C} \left\langle [\nu], [H(j)]  \right\rangle_\mathcal T [H(j)] \in V_p(\mathcal {T})  \end{equation}
and thus:
\begin{equation} \label{2}[\nu_p]=[\nu] \otimes [\nu] \otimes \cdots \otimes [\nu] =
\sum_{\sigma \in \mathcal S} 
( \prod_{i=1}^p \left\langle [\nu], [H(\sigma_i)]  \right\rangle_\mathcal T) [H(\sigma)] \in V_p(\mathcal {T}_p).\end{equation}

Again by  orthonormality, we have: 
\begin{equation} \label{3} [E]= \sum_{j \in C} \left\langle [E], [H(j)]  \right\rangle [H(j)] =\sum_{j \in C} \left\langle N'(j) \right\rangle_p [H(j)] \in V_p(\mathcal {T}) \end{equation}
\begin{equation} \label{4} [E_p]= \sum_{\sigma \in \mathcal{S}} \left\langle [E], [H(\sigma)]  \right\rangle [H(\sigma)] = \sum_{\sigma \in \mathcal{S}} \left\langle N'_p(\sigma) \right\rangle_p [H(\sigma)] \in V_p(\mathcal {T}_p),\end{equation}
 where $N'_p(\sigma)$ is $N'_p$ with the ith lift of the core colored $\sigma_i$ for all $i$. 
 Thus by \ref{1} and \ref{3}
  \begin{equation} \label{5} \left\langle N \right\rangle_p= \left\langle  [\nu], E \right\rangle_{\mathcal T}= \sum_{j \in C}   \left\langle [\nu], [H(j)]  \right\rangle_\mathcal T \left\langle N'(j) \right\rangle_p ,\end{equation}
and by \ref{2} and \ref{4}
   \begin{equation} \label{6} \left\langle N_p \right\rangle_p= \left\langle  [\nu_p], E_p \right\rangle_{\mathcal T_p}= \sum_{\sigma \in \mathcal{S}} ( \prod_{i=1}^p \left\langle [\nu], [H(\sigma_i)]  \right\rangle_\mathcal T) \left\langle N'_p(\sigma) \right\rangle_p .\end{equation}

 Note that 
 $N'_p(\sigma)$ is diffeomorphic to $N'_p(\tau (\sigma))$, and thus $\left\langle N'_p(\sigma) \right\rangle_p$ is constant on orbits of $\tau$. Also $ \prod_{i=1}^p \left\langle [\nu], [H(\sigma_i)]  \right\rangle_\mathcal T$ is constant on orbits of $\tau$.
 Since the non-constant orbits of $\tau$ have order $p$,  we have:  
 \begin{align}\label{7}  \left\langle N_p \right\rangle_p  &\equiv  \sum_{j \in C} ( \prod_{i=1}^p \left\langle [\nu], [H(j)]  \right\rangle_\mathcal T) \left\langle N'_p(\tilde j) \right\rangle_p\\ \notag
&\equiv  \sum_{j \in C}  \left( \left\langle [\nu], [H(\tilde j)]  \right\rangle_\mathcal T \right)^p \left\langle N'_p( j) \right\rangle_p
\pmod{ p \ \BZ[ \frac 1 {2r} ,\xi]} . \end{align}

As \cite[Theorem 1$'$]{G2}  is already proved for simple covers and $N'_p(j)$ is a simple cover of 
$N'(j)$, we have that:
 
\begin{equation}  \left\langle N'_p( j) \right\rangle_p \equiv  (\left\langle N'( j) \right\rangle_p)^p \pmod{ p \ \BZ[ \frac 1 {2r} ,\xi]} .\end{equation}

Substituting this in \ref{7} and comparing with \ref{5}, the result follows 
 \end{proof}

\end{document}